\documentclass[11pt]{article}
\usepackage{amssymb}
\usepackage{amsmath} 
\usepackage{setspace}
\begin{document}

\newtheorem{theorem}{Theorem}[section]
\newtheorem{tha}{Theorem}
\newtheorem{conjecture}[theorem]{Conjecture}
\newtheorem{corollary}[theorem]{Corollary}
\newtheorem{lemma}[theorem]{Lemma}
\newtheorem{claim}[theorem]{Claim}
\newtheorem{proposition}[theorem]{Proposition}
\newtheorem{construction}[theorem]{Construction}
\newtheorem{definition}[theorem]{Definition}
\newtheorem{question}[theorem]{Question}
\newtheorem{problem}[theorem]{Problem}
\newtheorem{remark}[theorem]{Remark}
\newtheorem{observation}[theorem]{Observation}

\newcommand{\ex}{{\mathrm{ex}}}

\def\endproofbox{\hskip 1.3em\hfill\rule{6pt}{6pt}}
\newenvironment{proof}%
{%
\noindent{\it Proof.}
}%
{%
 \quad\hfill\endproofbox\vspace*{2ex}
}
\def\qed{\hskip 1.3em\hfill\rule{6pt}{6pt}}

\def\ex{{\rm\bf ex}}
\def\cA{{\cal A}}
\def\cB{{\cal B}}
\def\cC{{\cal C}}
\def\cE{{\cal E}}
\def\cF{{\cal F}}
\def\cG{{\cal G}}
\def\cH{{\cal H}}
\def\ck{{\cal K}}
\def\cI{{\cal I}}
\def\cJ{{\cal J}}
\def\cL{{\cal L}}
\def\cM{{\cal M}}
\def\cN{{\cal N}}
\def\cP{{\cal P}}
\def\cQ{{\cal Q}}
\def\cS{{\cal S}}
\def\pkl{\mathcal{P}^{(k)}_\ell}
\def\podd{\mathcal{P}^{(k)}_{2t+1}}
\def\peven{\mathcal{P}^{(k)}_{2t+2}}
\def\cT{{\mathcal T}}

\voffset=-0.5in

\pagestyle{myheadings}
\markright{{\small \sc F\"uredi:}
  {\it\small Linear trees in uniform hypergraphs}}

\title{\huge\bf Linear trees in uniform hypergraphs}

\author{
Zolt\'an F\"uredi\thanks{R\'enyi Institute of Mathematics, POB 127, Budapest, 1364 Hungary.
E-mail: z-furedi@illinois.edu
\newline
Research supported in part by the Hungarian National Science Foundation
 OTKA 104343, 
 and by the European Research Council Advanced Investigators Grant 267195.
\newline ${}$
   \quad   
    {\rm\small {\jobname}.tex,}
   \hfill
Version as of
    March 21, 2012. Slightly revised on May 26, 2013.
}}

\date{}

\maketitle
\begin{abstract}
Given a tree $T$ on $v$ vertices 
 and an integer $k\geq 2$ one can define the $k$-expansion $T^{(k)}$  as a $k$-uniform linear
 hypergraph by enlarging each edge with a new, distinct set of $k-2$ vertices.
$T^{(k)}$ has  $v+ (v-1)(k-2)$ vertices.
The aim of this paper is to show that using the delta-system method one can easily
 determine  asymptotically the size of the largest $T^{(k)}$-free $n$-vertex hypergraph, i.e.,
 the Tur\'an number of $T^{(k)}$.
\end{abstract}


\section{Definitions: kernel-degree, Tur\'an number}

A {\it hypergraph} $H=(V,\cF)$ consists of a set $V$ of vertices and a
 set $\cF=E(H)$ of edges, where each edge is a subset of $V$.
We call the edges of $H$ {\it members} of $\cF$.
We say that $H$ is a {\it $k$-uniform hypergraph} or $\cF$ is a {\it $k$-uniform set system}
 if each member of $\cF$ is a $k$-subset of $V$.
To simplify notation we frequently identify the hypergraph $H$ to its edge set $\cF$.
If $|V|=n$, it is often convenient to just let $V=[n]=\{1,\ldots, n\}$.
We also write $\cF\subseteq {V\choose k}$ to indicate that $\cF$ is a $k$-uniform hypergraph,
 or $k$-{\it graph} for short, on vertex set $V$.
So ${V\choose k}$ 
  denotes the {\it complete} $k$-graph on vertex set $V$.
A set $S\subseteq V$ is a {\it transversal} (or vertex-cover) of the (hyper)graph $H=(V, {\mathcal E})$
 if $S\cap E\neq \emptyset$ for all $E\in {\mathcal E}$.
Let $\tau(H)$ denote the minimum number of vertices to cover all edges of  $H$, i.e., the
 {\it transversal} number of $H$.
A set of edges $\cM \subseteq E(H)$ is called a {\it matching} if it consists of disjoint members of $E(H)$.
$\nu(H)$ denotes the {\it matching number} of $H$, i.e., the maximum number of pairwise disjoint
 edges of $H$. 
A family of sets $\{F_1,\ldots, F_s\}$ is said to form a {\it $\Delta$-system} of size $s$
 with {\it kernel} $C$ if $F_i\cap F_j=C$ for all $1\leq i<j \leq s$.

Given a family $\cF\subseteq {[n]\choose k}$ and a subset $W\subseteq [n]$,
we define the {\it degree} of $W$ in $\cF$ as
$$\deg_\cF(W)=|\{F: F\in \cF, W\subseteq F\}|.$$
The hypergraph $\{F: F\in \cF, W\subseteq F\}$ is denoted by  $\cF[W]$, 
 so $\deg_\cF(W)= |\cF[W]|$ and $\deg_\cF(\emptyset)=|\cF|$.

We define the {\it kernel degree} of $W$, denoted by $\deg^*_\cF(W)$, as
$$\deg^*_\cF(W)=\max \{s: \exists \mbox{ a $\Delta$-system of size  with kernel $W$  in } \cF\}.$$
In other words,  $\deg^*_\cF(W)$ is the matching number of
 $\{ E\setminus W: W\subset E\in \cE \}$. 

Given a family $\cH=\{ H_1, H_2, \dots\}$ of hypergraphs,
 the {\it $k$-uniform hypergraph Tur\'an number} of $\cH$,
 denoted by $\ex(n,\cH)$,
 is the maximum number of edges in a $k$-uniform hypergraph $\cF$ on $n$ vertices
 that does not contain a member of $\cH$ as a subhypergraph.
If we want to emphasize $k$, then we write $\ex_k(n, \cH)$.
An $\cH$-free family $\cF\subseteq {[n]\choose k}$ is called {\it extremal}  if $|\cF|=\ex(n,\cH)$.
If $\cH$ consists of a single hypergraph $H$, we write $\ex(n,H)$ for $\ex(n,\{H\})$.
Surveys on Tur\'an problems of graphs and hypergraphs can be found in~\cite{ZF} and~\cite{keevash}.

It is easy to show (see, e.g., Bollob\'as~\cite{bollobas}, p.~xvii, formula (0.5))
 that any graph $G=(V,{\mathcal E})$ with more than $(\delta -1) |V|$ edges
 contains an induced subgraph $G'$ with minimum degree at least $\delta$.
Then $G'$ contains every tree of $\delta +1$ vertices. We have
\begin{equation}\label{eq:11}
  \ex(n,T)\leq (v-2)n,
  \end{equation}
  where $T$ is any $v$-vertex forest, $v\geq 2$.

For integers $b\geq a\geq 0$, $b\geq t\geq 1$ we have
$$
{a\choose t}=\frac{a}{t}{a-1\choose t-1}\leq \frac{a}{t}{b-1\choose t-1}=\frac{a}{b}{b\choose t}.
   $$
This implies the following lemma.

\begin{lemma}\label{le:63}
Suppose that $z_1\geq z_2\geq \dots \geq z_m$ and $t$ are non-negative integers, $z_1 \geq t\geq 1$.
Then
\begin{equation}\label{eq:61}
  \sum_{1\leq i\leq m} {z_i\choose t} \leq \frac{\sum z_i} {z_1} \binom{z_1}{t}.
  \end{equation}
  \end{lemma}

\section{Preliminaries: matchings, paths, stars}

The Erd\H{o}s-Ko-Rado theorem says that for $n\geq 2k$ the maximum size
 of a $k$-uniform family on $n$ vertices in which every two members intersect
 is ${n-1\choose k-1}$, with equality achieved by taking all the subsets of $[n]$ containing a fixed element.
If we let $M^{(k)}_\nu$ denote the $k$-uniform hypergraph consisting of $\nu$ disjoint $k$-sets,
 then the Erd\H{o}s-Ko-Rado theorem says $\ex_k(n,M^{(k)}_2)={n-1\choose k-1}$ for $n\geq 2k$.
More generally, Erd\H{o}s~\cite{erdos-matching} showed for any  positive integers $k,\nu$ 
 there exists a number $n(k,\nu)$ such that the following holds.
For all $n>n(k,\nu)$, if $\cF\subseteq {[n]\choose k}$ contains no $\nu+1$ pairwise disjoint
 members then
\begin{equation}\label{eq:21}
 |\cF|\leq {n\choose k}-{n-\nu\choose k}.
  \end{equation}
Furthermore, the only extremal family $\cF$ consists of all the $k$-sets of $[n]$
 meeting some fixed set $S$ of $\nu$ elements of $[n]$. 

The value of $n(2,\nu)$ was determined by Erd\H{o}s and Gallai~\cite{erdos-gallai}. 
Frankl, R\"odl, and Rucin\'ski~\cite{FRR} showed $n(3,\nu)\leq 4\nu$.
Finally,  $n(3,\nu)$ was determined by {\L}uczak and Mieczkowska~\cite{LM} for 
 large $\nu$  (for $\nu> 10^5$), and by Frankl~\cite{F12} for all $\nu$.
In general, Huang, Loh, and Sudakov~\cite{HLS} showed $n(k,\nu)< 3\nu k^2$,
 which was slightly improved in~\cite{FLM} and greatly improved to $n(k,\nu)\leq (2\nu+1)k-\nu$
 by Frankl~\cite{frankl-matching}. 
Summarizing, for fixed $k$ and $\nu$ as $n\to \infty$ we have that
\begin{equation}\label{eq:23}
  \ex_k(n,M^{(k)}_\nu)=(\nu +o(1)){n-1 \choose k-1}.
  \end{equation}

A  {\it linear path}  of length $\ell$ is a family of sets $\{F_1,\ldots, F_\ell\}$
 such that $|F_i\cap F_{i+1}|=1$ for each $i$ and $F_i\cap F_j=\emptyset$ whenever $|i-j|>1$.
Let $\pkl$ denote the $k$-uniform linear path of length $\ell$.  It is unique up to isomorphisms.
Note that this notation is different from what is usually used, where $P_v$ denotes a $v$-vertex path.
Concerning the graph case ($k=2$) Erd\H{o}s and Gallai~\cite{erdos-gallai} proved that
 $\ex_2(n, \cP_\ell^{(2)}) \leq \frac{1}{2}(\ell-1)n$. Here equality holds if $G$ is the disjoint
  union of complete graphs on $\ell$ vertices.
The value of $\ex_2(n, \cP_\ell^{(2)})$
 was determined for all $n$ by Woodall~\cite{Woo} and Kopylov~\cite{Kop}. 

Concerning linear paths of two edges
Erd\H os and S\'os~\cite{ESos} proved for triple systems ($k=3$) that $\ex_3(n, \cP_2^{(3)})=n$ or $n-1$
 (according to $n$ is divisible by 4 or not and $n\geq 4$).
They conjectured that
\begin{equation}\label{eq:24}
  \ex_k(n, \cP_2^{(k)})={n-2\choose k-2}
   \end{equation}
for $k\geq 4$ and
 sufficiently large $n$ with respect to $k$, and this was proved by Frankl~\cite{frankl1}.
The case $k=4$ was finished for all $n$ by Keevash, Mubayi, and Wilson~\cite{KMW}.

The case $\ell<k$ was asymptotically determined in~\cite{frankl-furedi}.

Since the paper of G. Y. Katona and Kierstead~\cite{KK} (1999) there is a renewed interest
 concerning paths and (Hamilton) cycles in uniform hypergraphs.
Most of these are Dirac type results
 (large minimum degree implies the existence of the desired substructure)
 like in K\"uhn and Osthus~\cite{KO},  R\"odl, Ruci\'nski,  and Szemer\'edi~\cite{RRS}.

The present author, Tao Jiang, and Robert Seiver~\cite{FJS} determined $\ex_k(n, \pkl)$ exactly,
for {\bf all} fixed $k,\ell$, where $k\geq 4$, and sufficiently large $n$ proving
\begin{equation}\label{eq:25}
\ex_k(n,\podd)={n-1\choose k-1}+{n-2\choose k-1}+\ldots+{n-t\choose k-1},
  \end{equation}
where the only extremal family consists of all
the $k$-sets in $[n]$ that meet some fixed set $S$ of $t$ elements, and 
\begin{equation}\label{eq:26}
\ex(n,\peven)={n-1\choose k-1}+{n-2\choose k-1}+\ldots+{n-t\choose k-1}+{n-t-2\choose k-2},
   \end{equation}
 where the only extremal family consists of all
the $k$-sets in $[n]$ that meet some fixed set $S$ of $t$ elements plus all the $k$-sets in $[n]\setminus S$
that contain some two fixed elements. 
`Sufficiently large' $n$ means that \eqref{eq:25} and \eqref{eq:26} hold when $kt=O(\log\log n)$.
It is \textbf{conjectured} that they hold for all (or at least almost all) $n$'s.
The method in~\cite{FJS} does not quite work for  the $k=3$ case (cf. the remark after Lemma~\ref{le:62} below)
but it is \textbf{conjectured} that still a similar result holds for $k=3$.

A ({\it linear}) {\it star} of size $\ell$ with center $x$ is a family of sets $\{F_1,\ldots, F_\ell\}$
 such that  $x \in F_i$ for all $i$ but the sets   $F_i \setminus \{ x\}$ are pairwise disjoint.
Let $\cS_\ell^{(k)}$ denote the $k$-uniform star of size $\ell$.
It is obvious that $\ex_2 (n,\cS_\ell^{(2)})= \lfloor (\ell-1)n/2 \rfloor$ (for $n \geq \ell$).
Chung and Frankl~\cite{CF} gave an exact formula for $\ex_3 (n,\cS_\ell^{(3)})$ for $n> 3\ell^3$.
The following asymptotic was proved for any fixed $\ell\geq 2$, $k\geq 5$ in~\cite{frankl-furedi}.
\begin{equation}\label{eq:27}
  \ex_k(n,\cS_\ell^{(k)})=(\varphi(\ell) +o(1)){n-2 \choose k-2},
  \end{equation}
where $\varphi(\ell)=\ell^2-\ell$ for $\ell$ is odd and
 it is $\ell^2-\frac{3}{2}\ell$ when $\ell$ is even.
According to the above mentioned result of Chung and Frankl \eqref{eq:27} holds for $k=3$ too.
The order of magnitude $\ex_4(n,\cS_\ell^{(4)})=\Omega(\ell^2n^2)$ was also proven in~\cite{frankl-furedi},
 and it is \textbf{conjectured} that \eqref{eq:27} holds for $k=4$ too.

\section{Generalized $k$-forests, an upper bound}

Let us define a generalized $k$-{\it forest} in the following inductive way.
Every $k$-graph consisting of a single edge is a $k$-forest.
Suppose that $\cT= \{ E_1, E_2, \dots, E_u\} \subseteq {V\choose k}$ is a $k$-forest
 and suppose that $A:=A_{u+1}\subset E_i$ for some $1\leq i\leq u$,
 and  $B\cap V=\emptyset$, $|A|+|B|=k$, then $\{ E_1, E_2, \dots, E_u, E_{u+1} \}$
 is a $k$-forest with $E_{u+1}:=A\cup B$.
If it is connected then it is called a generalized $k$-tree.
In that case all defining sets $A_2, \dots, A_{u+1}$ are nonempty.
For graphs ($k=2$) the above process leads to the usual notions of
 forests and trees.
If each defining set $A_i$ is a singleton or empty then we obtain a {\it linear} forest,
if each defining set is either empty or has $k-1$ elements, then we get a {\it tight}
 forest.
A forest $\cT$ of $q$ edges has at least $q+k-1$ vertices and here equality holds
 if and only if $\cT$ is a tight $k$-tree.

Consider a $k$-forest  $\cT= \{ E_1, E_2, \dots, E_q\}$.  
If a defining set $A_{u+1}\subset E_i$ for some $1\leq i\leq u< q$ is smaller than $k-1$, 
  then take an
 element $x\in (E_i\setminus A_{u+1})$ and another one $y\in (E_{u+1}\setminus A_{u+1})$ 
 and place the new $k$-set $E:= E_i\setminus \{ x\} \cup \{ y\}$ between $E_u$ and $E_{u+1}$. 
The new sequence of $k$-sets $\{E_1, \dots, E_{u}, E, E_{u+1}, \dots, E_q\}$
 is again a $k$-forest with the same defining sets except we add $E_i\setminus \{ x\}=E\setminus \{ y\}$ 
 to the list for $E$ and replace $A_{u+1}$ by $(A_{u+1}\cup \{ y\})$
 and use the relation $(A_{u+1}\cup \{ y\}) \subset E$ for $E_{u+1}$. 
Repeating this process we obtain the following statement.

\begin{proposition}\label{prop:31}
Suppose that $\cT$ is a generalized $k$-forest of $v$ vertices. Then
 there is a tight $k$-tree $\cT^+$ on the same vertex set such that 
 $\cT$ is a subfamily of $\cT^+$. 
  \end{proposition}

We are going to prove the following upper bound for the Tur\'an number of $k$-forests. 

\begin{theorem}\label{th:31}
Suppose that $\cT$ is a generalized $k$-forest of $v$ vertices. Then
\begin{equation}\label{eq:31}
    \ex_k(n, \cT)\leq (v-k){n \choose k-1}.
     \end{equation}
  \end{theorem}

\noindent
\emph{Proof.}\enskip
By the previous Proposition, it is enough to prove the case when 
 $\cT$ is a tight $k$-forest.
 
Suppose that $\cH\subseteq {[n]\choose k}$ avoids the tight $k$-forest $\cT=
\{ E_1, \dots, E_q \}$, we have $q=v-k+1$. 
Set $A_{i}:=E_{i}\cap (E_1\cup \dots\cup E_{i-1})$,
 $2\leq i\leq q$.
We have that $A_{i}\subset E_{\alpha(i)}$ for some $1\leq \alpha(i)< i$, $|A_i|=k-1$. 
Define a list of hypergraphs $\cH_0:=\cH \supset \cH_1\supset \dots \supset \cH_m$
 and sets $X_1, \dots, X_m$, 
 as follows.

If $\cH_m=\emptyset$ we stop.
If one can find a set $X\subset [n]$ such that $|X|=k-1$ and $\deg_{\cH_m}(X)\leq (v-k)$
 then let $X_{m+1}:=X$ and $\cH_{m+1}:= \cH_m \setminus \cH_m[X]$.
If there is no such set $X$ then we stop.

We claim that $\cH_m$ should be the empty family.
Otherwise, we can embed $\cT$ into $\cH_m$ as follows.
Start with any edge $E_1\in \cH_m$.
We define the other edges $E_2, \dots, E_q$ one by one.
Observe that for any $(k-1)$-element subset $X$, $X\subsetneq E\in \cH_m$ we have $\deg_{\cH_m}(X)\geq v-k+1$.
Suppose that $E_1, \dots, E_u$ had already been defined together with $A_2, \dots, A_u$, and $u<q$.
Locate $A_{u+1}$ in $E_1\cup \dots \cup E_u$.
Since $\deg_{\cH_m}(A_{u+1})\geq (v-k+1) > |E_1\cup \dots \cup E_u|-|A_{u+1}|$ there is an 
 $E:=E_{u+1}\in \cH_m[A_{u+1}]$ such that $E\setminus A_{u+1}$ is disjoint to $E_1\cup \dots \cup E_u$. 

In the sequence $X_1, \dots, X_m$ there is no repetition, so we get
\begin{equation*}
  |\cH| =\sum_{i} \deg_{\cH_i}(X_i) \leq  (v-k) {n \choose k-1}. \quad \square
  \end{equation*}

Note that Theorem~\ref{th:31} gives the correct order of magnitude if
$\cap \cT=\emptyset$, since then ${n-1\choose k-1}$ is a lower bound.
However, the determination of the best coefficient of the binomial term seems to be
 extremely difficult.
Erd\H os and S\'os \textbf{conjectured} for graphs (i.e., $k=2$) and
 Kalai 1984 for all $k$, see in~\cite{frankl-furedi}, that for a $v$-vertex tight tree $\cT$
$$ \ex_k(n,\cT)\leq \frac{v-k}{k}{n \choose k-1}.
  $$
For any given tight tree $\cT$ a matching lower bound, i.e.,  $(1-o(1))$ times the conjectured upper bound,
 can be given for $n \to \infty$ as follows.
Consider a $P(n,v-1,k-1)$ packing $P_1$, \dots, $P_m$ on the vertex set $[n]$
  (i.e., $|P_i|=v-1$ and $|P_i\cap P_j|< k-1$ for $1\leq i<j\leq m$) and replace
  each $P_i$ by a complete $k$-graph.
We obtain a $\cT$-free hypergraph.
Then  R\"odl's~\cite{R} theorem on almost optimal packings gives 
$$
   \ex_k(n,\cT)\geq (1-o(1))\frac{{n\choose k-1}}{{v-1\choose k-1}}\times
  {v-1\choose k} = (1+o(1))\frac{v-k}{k}{n \choose k-1}.
  $$ 

The Erd\H os--S\'os conjecture has been recently proved by a monumental work of 
 Ajtai, Koml\'os, Simonovits, and Szemer\'edi~\cite{AKSSzSharp}, for $v\geq v_0$.

The Kalai conjecture has been proved for {\it star-shaped} $k$-trees
 in~\cite{frankl-furedi}, i.e., whenever $\cT$ contains a central edge which 
 intersects all other edges in $k-1$ vertices. 
For $k=2$ these are the diameter 3 trees,  'double stars'.

There is only one more class of $k$-trees where the exact asymptotic is know,
 namely what is called an {\it intersection condensed family}.
For such a $\cT$ we denote $|\cap \cT|$ by $p_\infty$, and the number of vertices
 of degree at least two by $p_2$ and suppose that $2p_\infty+p_2+2\leq k$
  (Theorem~5.3 in~\cite{frankl-furedi}).

There are many different definitions of a `path' in a hypergraph.
Gy\H ori, G. Y. Katona, and Lemons~\cite{gyori} determined the exact value
 of the Tur\'an number of the so-called
  {\it Berge}-paths for infinitely many $n$'s.
Mubayi and Verstra\"ete~\cite{MV1} gave good bounds for the Tur\'an number of
 $k$-uniform {\it loose} paths of length $\ell$.

The aim of this paper is to present the best coefficient for
 a wide class of linear trees, thus generalizing
 the results in the previous section about matchings~\eqref{eq:23},
 paths~\eqref{eq:24}--\eqref{eq:26} and stars~\eqref{eq:27}.

\section{The main result, finding expanded forests in $k$-graphs}

Given a graph $H$, the {\it $k$-blowup} (or $k$-expansion),
 denoted by $[H]^{(k)}$   (or $H^{(k)}$ for short),
is the $k$-uniform hypergraph obtained
from $H$ by replacing each edge $xy$ in $H$ with a $k$-set $E_{xy}$ that consists of $x,y$ and $k-2$ new
vertices such that for distinct edges $xy, x'y'$, $(E_{xy}-\{x,y\})\cap (E_{x'y'}-\{x',y'\})=\emptyset$.
If $H$ has $p$ vertices and $q$ edges, then $H^{(k)}$ has $p+q(k-2)$ vertices and $q$ hyperedges.
The resulting $H^{(k)}$ is a $k$-uniform 
 hypergraph whose vertex set contains the vertex set of $H$. 

Given a forest $T$ define the following
\begin{equation}\label{eq:41}
\sigma (T):=\min \{  |X|+ e(T\setminus X): X\subset V(T) \textrm{ is independent in } T\}.
 \end{equation}
Here $T\setminus X$ is the forest left from $T$ after deleting the vertices of $X$ and the
 edges incident to them, $e(G)$ stands for the number of edges of the graph $G$. 
Since the edges avoiding $X$ can be covered one by one we have that $\tau(T)\leq \sigma(T)$
 but here equality should not hold.
For example, if $T$ consists of a path of four vertices $a_1b_1b_2a_2$ with
 $2d+2c$ pendant edges such that $d>c\geq 1$ and each $a_i$ has $d$ degree-one neighbors
  and each $b_i$ has $c$ of those, then one can easily see that $\tau(T)=4$ but $\sigma(T)=2c+3$.

\begin{theorem}\label{th:41}
Given a forest  $T$ with at least one edge and an integer $k\geq 4$. Then we have as $n\to \infty$, that
\begin{equation}\label{eq:42}
  \ex(n, T^{(k)}) = \left( \sigma(T) -1+o(1)\right){n \choose k-1}.
  \end{equation}
  \end{theorem}

Our result, naturally, gives the same asymptotic as
 Theorem~5.3 in~\cite{frankl-furedi} whenever both can be applied to $T^{(k)}$.
We \textbf{conjecture} that \eqref{eq:42} holds for $k=3$, too.

According to \eqref{eq:27} (and the remark after that) the above asymptotic 
 holds for stars, since the answer in this case is $o(n^{k-1})$.
For every other forest $\sigma\geq \tau\geq 2$. 

Let us note that Mubayi~\cite{dhruv} and Pikhurko~\cite{oleg} determined precisely
 (for large $n$)
  the Tur\'an number of the $k$-expansion of some other graphs, namely for the
   complete graph $K_\ell$ for $\ell> k\geq 3$.
For smaller values of $\ell$ we know that
  $\ex_k (n, K_3^{(k)})={n-1\choose k-1}$ for $n>n_0(k)$, $k\geq 3$, a former
   conjecture of Chv\'atal and Erd\H os, established in~\cite{frankl-furedi}. 
A few more related exact results can be found in~\cite{lale-furedi}.

\section{The product construction}\label{S:5}

Given two set systems (or hypergraphs) $\cA$ and $\cB$ their {\it join} 
 is the family $\{ A\cup B: A\in \cA, B\in \cB \} $. We denote this new hypergraph by
$\cA \Join \cB$.

Call a set $Y$ 1-{\it cross-cut}  of a family $\cC$ if $|Y\cap E|=1$ holds for each $E\in \cC$.
Define $\tau_1(\cC)$ as the minimum size of a 1-{\it cross-cut} of $\cC$ (if such
 cross-cut exists, otherwise $\tau_1:=\infty$).
We claim that for every forest $T$ and $k\geq 3$ the following holds.
\begin{equation}\label{eq:52}
  \sigma(T) =\tau_1(T^{(k)}). 
  \end{equation}
Indeed, suppose that $X \subset V(T)$ yields the minimum in \eqref{eq:41}.
Then $X$ is an independent set of $T^{(k)}$ avoiding $e(T\setminus X)$ edges of it.
Taking an element $x(E)\in (E\setminus V(T))$ from each such edge 
 and joining them to $X$ one gets a $1$-cross-cut of size
 $\sigma(T)$. We obtain $\tau_1\leq \sigma$.
On the other hand, if  $S$ is a $1$-cross-cut of $T^{(k)}$ and $|S|=\tau_1$, 
 then $X:= S\cap V(T)$ is an independent set in $T$
 and it avoids exactly $|S|-|X|$ edges, so $\sigma \leq |S|=\tau_1$. 
  
Thus $\sigma(T)$ is the minimum size of a set $Y$ such that
 $T^{(k)}$ can be embedded into ${Y \choose 1}\Join{Z\choose k-1}$
where $Y$ and $Z$ are disjoint sets.  
This means that in case of $Y:=[\sigma-1]$,  $Z:=[n]\setminus Y$ the hypergraph
${Y \choose 1}\Join{Z\choose k-1}$ does not contain any copy of $T^{(k)}$.
We obtain the lower bound
\begin{equation}\label{eq:53}
  \ex(n,T^{(k)})\geq |{Y \choose 1}\Join{Z\choose k-1}|=
    |\{ E: E\in {[n]\choose k},  |E\cap [\sigma-1]|=1\}|= (\sigma-1){n -\sigma+1 \choose k-1}.
\end{equation}

\section{The graph of 2-kernels, starting the proof with the delta-system method}

Given a family $\cF\subseteq {[n]\choose k}$,
  the {\it kernel-graph} with {\it threshold } $s$ is a graph $G:=G_{2,s}(\cF)$ on $[n]$
  such that $\forall x,y\in [n]$,  $xy\in E(G)$ if and only if $\deg^*_{\cF}(\{x,y\})\geq s$.
The following (easy) lemma shows the importance of this definition.

\begin{lemma}[see~\cite{FJS}]\label{le:61}
Let $H$ be a graph with $q$ edges, $s=kq$,
 and let $\cF\subseteq {[n]\choose k}$.
Let $G_2$ be the kernel graph of $\cF$ with threshold $s$.
If $H\subseteq G_2$, then $\cF$ contains a copy of $H^{(k)}$. 
\hfill $\square$
\end{lemma}

The {\em delta-system method},  started by  Deza, Erd\H os and Frankl~\cite{DEF}, 
  is a powerful tool for solving set system problems.
Using a structural lemma from~\cite{furedi-1983} and the method developed
 in~\cite{ff40,frankl-furedi} the following theorem was obtained in~\cite{FJS}
 (see Theorem~3.8 and the proof of Lemma~4.3 there).

\begin{lemma}[see~\cite{FJS}]\label{le:62}
Let $\cF\subseteq {[n]\choose k}$, $T$ a forest of $v$ vertices, $s=kv$, $G_2:=G_{2,s}(\cF)$,
 and suppose that $\cF$ does not contain $T^{(k)}$.
 Then there is a constant $c:=c(k,v)$ and a partition $\cF=\cF_1\cup \cF_2$ with the following properties.
 \newline --- \enskip $|\cF_1|\leq c {n-2\choose k-2}$.
 \newline --- \enskip Every edge $F\in \cF_2$ has a center (not necessarily unique) $x(F)\in F$
   such that $G_2|F$ contains a star of size $k-1$ with center $x(F)$.
   In other words, $\{x(F),y\}\in E(G_2)$ for all $y\in F\setminus \{ x(F)\}$.
\hfill $\square$
  \end{lemma}

Actually, the delta-system method describes the intersection structure of $\cF$
 in a more detailed way, but for our purpose this lemma will be sufficient.
The above lemma (and in fact the main result of this paper, Theorem~\ref{th:41})
 preceded \eqref{eq:25}--\eqref{eq:26},
 see~\cite{furedi_2011}, but since the proof of Lemma~\ref{le:62} is now available in~\cite{FJS}
  we omit the details here.

Note that this is the only point where $k\geq 4$ is used.
Lemma~\ref{le:62} is not true for $k=3$.
The 3-graph $\cF^3$ obtained by joining a matching of size $t$ and $t$ one-element sets has
  $n=3t$ vertices, $t^2= n^2/9= \Omega(n^{k-1})$ edges, it does not contain
  any linear tree except stars but $G_{2,s}(\cF^3)$ forms a matching for every $s\geq 2$.

\section{Proof of the Main Theorem}

Suppose that $\cF\subseteq {[n]\choose k}$ avoids the $k$-expansion
 of the $v$-vertex forest $T$, $k\geq 4$.
We are going to give an upper bound for $|\cF|$.
As noted above we may suppose that $T$ is not a star, $\sigma(T)\geq \tau(T)\geq 2$.

Define $s=vk$ and let $G_2$ be the kernel graph with threshold $s$ with respect
 to family $\cF$ as defined in the previous Section.
This graph avoids $T$ by Lemma~\ref{le:61}, so \eqref{eq:11} implies
\begin{equation}\label{eq:71}
  e(G_2)\leq (v-2)n.
  \end{equation}

Consider the degree sequence of $G_2$ and suppose that
  $$ \deg(x_1)\geq \deg(x_2)\geq \dots \geq \deg(x_{n-1})\geq \deg(x_n).
    $$
Let $L:= \{ x_1, \dots, x_\ell\}$ be the set of highest degrees.
We will define $\ell$ later as $n^\varepsilon$ so keep in mind that
 it is relatively large.
Using \eqref{eq:71} we obtain
\begin{equation}\label{eq:72}
    z:= \deg_{G_2}(x_{\ell+1})\leq \frac{\deg(x_1)+ \dots +\deg(x_{\ell+1})}{\ell+1}
      \leq \frac{2e(G_2)}{\ell+1} < \frac{2(v-2)n}{\ell}.
  \end{equation}

Consider the partition $\cF=\cF_1\cup \cF_2$ given by Lemma~\ref{le:62}.
Let $\cF_3$ be the edges of $\cF_2$ with center outside $L$.
Using Lemma~\ref{le:62} and~\eqref{eq:61} then \eqref{eq:71}, a triviality and \eqref{eq:72} we get
\begin{eqnarray} \notag
  |\cF_3|&\leq& \sum_{\ell+1 \leq i\leq n}  {\deg(x_i)\choose k-1}
    \leq \frac{\sum_i \deg(x_i)}{z} {z\choose k-1} \\
   {} &\leq &\frac{2(v-2)n}{z}{z \choose k-1}
     < \frac{2(v-2)n}{(k-1)!}z^{k-2}\leq \frac{2^{k-1}(v-2)^{k-1}}{(k-1)!}\frac{n^{k-1}}{\ell^{k-2}} .
\label{eq:73}
  \end{eqnarray}

Every edge of $\cF\setminus (\cF_1\cup \cF_3)$ meets $L$.
Let $\cF_4$ be the set of members of $\cF$ meeting $L$ in at least two vertices.
Obviously
\begin{equation}\label{eq:74}
  |\cF_4|\leq {\ell\choose 2}{n-2\choose k-2} \leq \frac{1}{2\times (k-2)!} \ell^2 n^{k-2}.
  \end{equation}

The edges of $\cF\setminus (\cF_1\cup \cF_3\cup \cF_4)$ meet $L$ in exactly one element.
Let $\cF_5$ be the family of edges of $\cF$ satisfying $|F\cap L|=1$ and $\deg_\cF(F\setminus L)\leq \sigma-1$.
Obviously,
\begin{equation}\label{eq:75}
  |\cF_5|\leq (\sigma-1){n-\ell\choose k-1}.
  \end{equation}

The rest of the edges, i.e., those from $\cF_6:=\cF\setminus (\cF_1\cup \cF_3\cup \cF_4\cup \cF_5)$
 are of the form $F=\{ a\}\cup B$ where $a\in L$, $B\cap L=\emptyset$
 and $\deg_\cF(F\setminus L)\geq \sigma$.
For every set $A\in {L\choose \sigma}$ define
 $\cB_A$ as the  $k-1$ uniform family
 $$
   \cB_A:=\{ B: 
     \{ a\} \cup B\in \cF\mbox{ for all }a\in A\}.
   $$
Also set
$$
 \cF_A:=\{ F\in \cF: a\in A, B\in \cB_A,\mbox{ and }\{ a\} \cup B=F\}.
 $$
We have $\cF_6\subseteq \cup_A \cF_A$ where $|A|=\sigma$, $A\subseteq L$.

Consider $T^{(k)}$.
As noted in Section~\ref{S:5}, there is a 1-cross-cut, a set $Y$ of size $\sigma$ meeting
 each $k$-edge of $T^{(k)}$ in a singleton.
Let $\cC$ be the $(k-1)$-uniform hypergraph obtained by deleting the elements of $Y$
 from the edges of $T^{(k)}$, $\cC:=\{E\setminus Y: E\in E(T^{(k)}) \}$.
Since $\cF_A$ does not contain $T^{(k)}$ we have that
  $\cB_A$ can not contain $\cC$ as a subhypergraph.
Also, $\cC$ is a generalized forest of at most $v-1$ edges so Theorem~\ref{th:31} gives
 $|\cB_A|\leq (v-2)(k-1){n\choose k-2}$.
We obtain
\begin{equation}\label{eq:76}
 |\cF_6|\leq \sum_{A\in {L \choose \sigma}} |\cF_A|
   = \sigma \sum_{A\in {L \choose \sigma}} |\cB_A|\leq \sigma {\ell \choose \sigma}
        (v-2)(k-1){n\choose k-2}.
  \end{equation}

Finally, since $\cF_2\subseteq \cF_3\cup \cF_4\cup \cF_5\cup \cF_6$ we have
$$
  |\cF|\leq |\cF_1|+ |\cF_3|+|\cF_4|+|\cF_5|+|\cF_6|.
   $$
Using the first part of Lemma~\ref{le:62}, \eqref{eq:73}, \eqref{eq:74}, \eqref{eq:75}
 and \eqref{eq:76} we obtain
\begin{equation}\label{eq:77}
  |\cF|\leq O(n^{k-2})+ O(\frac{n^{k-1}}{\ell^{k-2}})+ O(\ell^2 n^{k-2})+
 (\sigma-1){n-\ell\choose k-1} +O(\ell^\sigma n^{k-2}).
  \end{equation}
Defining $\ell\sim n^{1/(\sigma+1)}$ we obtain that
the sum of the $O()$ terms in \eqref{eq:77} is $O(n^{(k-1)-1/(\sigma+1)})=o(n^{k-1})$ and we are done.

\section{Further problems}

With a refined version of the above proof one can see that
  $$
  \ex(n, T^{(k)})= (\sigma-1){n \choose k-1}+O(n^{k-2}).
  $$
It seems to be a solvable problem to determine the
 exact value of this Tur\'an number (for $n> n_0(T,k)$) as it was done 
  for linear paths (for $k\ge 4$) in~\cite{FJS}, and 
 for linear cycles (for $k\ge 5$ only) in~\cite{i227}.
The forthcoming manuscript~\cite{FJ_partial_forest}
 generalizes these to a class of expanded forests, but most of the 
cases remain unsolved.


\end{document}